\documentclass{article}

\usepackage{arxiv}

\usepackage[utf8]{inputenc} 
\usepackage[T1]{fontenc}    
\usepackage{hyperref}       
\usepackage{url}            
\usepackage{booktabs}       
\usepackage{amsfonts}       
\usepackage{nicefrac}       

\title{On irrational values of the error function and gamma function }

\author{Ali Chtatbi}

\begin{document}
\maketitle

\begin{abstract}
In this paper, we establish the irrationality of some open problems in mathematics based on using a recursive formula that generate the complete sequence of numbers. see \cite{fridman2019prime}
But before getting into that we begin with some Ramanujan notable work of infinite series and continued fraction.
\end{abstract}

\keywords{Irrational numbers  \and Error function \and Gamma function }

\section{Introduction}
Irrational numbers are real numbers that cannot be constructed from ratios of
integers. Among the set of irrational numbers, two famous constants are $e$ and $\pi $. Lambert was the first mathematician that showed the irrationality of $\pi$ using continued fraction \cite{wallisser2011lambert}, years after Euler prove that e has non-repeating continued fraction expansion implies e is irrational \cite{euler1744fractionibus}, In 1815 Joseph Fourier introduces
a nice proof depends on the series representation of e \cite{de1815melanges}. Therefore, to reach
a contradiction that determines the truth of the proposition be false, but proof of
irrationality is still difficult to obtain. For instance ,the proofs of irrationality of $\pi+e$ or $\pi-e$ are open problems .

\section{Ramanujan remarkable formula :}
Ramanujan found the following remarkable formula proven in \cite{berndt2012ramanujan} and \cite{kouba2006inequalities}  :
\[\sqrt{\frac{\pi e^x}{2x}}=\frac{1}{x+}\frac{1}{1+}\frac{2}{x+}\frac{3}{1+}\frac{4}{x+}...+\left\{1+\frac{x}{1\cdot3}+\frac{x^2}{1\cdot3\cdot5}+\frac{x^3}{1\cdot3\cdot5\cdot7}+...\right\}.\]
This formula above is one of the many staggering relationships found by mathematical genius Ramanujan it has been described by some mathematicians as Ramanujan's most beautiful formula as it shows a remarkable connection between an infinite series, a continued fraction, and two of the fundamental mathematical constants e and $\pi$.
The second series of Ramanujan formula seems to be familliar:
\[1+\frac{x}{3}+\frac{x^2}{3\cdot5}+\frac{x^3}{3\cdot5\cdot7}+...=\sum_{n=0}^\infty \frac{x^n}{(2n+1)!!}\]
where $n!!$ is the double factoriel. by putting $x=-1$ we get 
\[\sum_{n=0}^\infty \frac{1}{(2n+1)!!}=\sqrt{\frac{e\pi}{2}}{erf}\left(\frac{1}{\sqrt 2}\right)\ \]
\[{{erf}(z)=\frac{2}{\sqrt{\pi}}\int_0^z e^{-t^2}dt}  \]
 where erf(z) is the error function                  The first term in Ramanujan's formula is the continued fraction is given by
\[\frac{1}{1+}\frac{1}{1+}\frac{2}{1+}\frac{3}{1+}\frac{4}{1+}...=
\sqrt{\frac{e \pi}{2}} erfc(\left(\frac{1}{\sqrt{2}}\right)=\sqrt{\frac{e \pi}{2}}\left(1-erf(\left(\frac{1}{\sqrt{2}}\right)\right) \] 
Where erfc(z) is the Complementary error function, defined as \[{{erfc}(z)=\frac{2}{\sqrt{\pi}}\int_0^z e^{-t^2}dt=1-{{erf}(z)}} \]
\subsection{Theorem : } For every integer $ \alpha > 1 , \beta> 0     $ the number  :
\[ \sum_{k=0}^{\infty} \frac{(\alpha k+\beta)-1}{\prod_{k=0}^{k-1}(\alpha k+\beta)}=  \\
\left(\frac{1}{\alpha}\right)^{2-\beta / \alpha}\left(\left(\frac{1}{\alpha}\right)^{\beta / \alpha-2}+\right.  
\sqrt[\alpha]{e} \alpha^{2} \Gamma\left(\frac{\beta}{\alpha}\right)+ \\
\sqrt[\alpha]{e} \alpha \Gamma\left(\frac{\beta}{\alpha}\right)-
\sqrt[\alpha]{e} \alpha \beta \Gamma\left(\frac{\beta}{\alpha}\right)+ 
\sqrt[\alpha]{e} \beta^{2} 
 \Gamma\left(\frac{\beta-\alpha}{\alpha}\right)+ \] \[ 
\sqrt[\alpha]{e} \alpha \Gamma\left(\frac{\beta-\alpha}{\alpha}\right)- \\
\sqrt[\alpha]{e} \beta \Gamma\left(\frac{\beta-\alpha}{\alpha}\right)- \\
\sqrt[\alpha]{e} \alpha \beta \\  \Gamma\left(\frac{\beta-\alpha}{\alpha}\right)- \\
\sqrt[\alpha]{e} \alpha^{2} \Gamma\left(\frac{\beta}{\alpha}, \frac{1}{\alpha}\right)- \\ 
\sqrt[\alpha]{e} \alpha \Gamma\left(\frac{\beta}{\alpha}, \frac{1}{\alpha}\right)+ \\  
\sqrt[\alpha]{e} \alpha \beta \\ \] \[ 
\Gamma\left(\frac{\beta}{\alpha},  \frac{1}{\alpha}\right)- \\
\sqrt[\alpha]{e} \beta^{2} \\
\Gamma\left(\frac{\beta}{\alpha}-1, \frac{1}{\alpha}\right)- \\
\sqrt[\alpha]{e} \alpha \Gamma\left(\frac{\beta}{\alpha}-1,\right. \\ 
\left.\quad \frac{1}{\alpha}\right)+\sqrt[\alpha]{e} \beta \\
 \Gamma\left(\frac{\beta}{\alpha}-1, \frac{1}{\alpha}\right)+ \\ \sqrt[\alpha]{e} \alpha  \beta\left.\Gamma\left(\frac{\beta}{\alpha}-1,\frac{1}{\alpha}\right)\right) \notin \mathbf{Q} \]
We gonna Proof the theorem by using Bertrand's postulate  states that if $k$>$3$ there is always at least one prime $ p_{k} $ between $k$ and $ 2k-1$. Equivalently, $ p_{k-1}< p_{k} < 2p_{k-1} $ . we can deduce that  $ (\alpha k + \beta) < f_{k} < (\alpha k + \beta) + \alpha $ for all k. we write $ f_{k} = (\alpha k + \beta) + r_{k} $, where $0< r_{k} < 1 $
We assume $ A(\alpha , \beta ) $ is rational. Then, there exist some   $ a , b \in N $ with $ a > 1 $ , such that  \[  A(\alpha, \beta )   =\frac{a}{b}\]    using the recurrance relation see [1] 
\[ f_{k+1} = (\alpha k + \beta )  \left (f_{k} - ( \alpha k +\beta )  \right ) = (\alpha k + \beta ) \left(r_{k}+1 \right ) \]  The statement give us as follows.  We multiply the above equation  by b. Notice that   $ b f_{k} $ is an integer for all k . in particular $ r_{k} >  1/b $ for all k .
By rearranging the above expression we get  : 
\[ \left(r_{k}+1 \right ) = \frac{f_{k+1}}{(\alpha k + \beta )} = \frac{(\alpha k+ \beta)+\alpha )}{(\alpha k + \beta )} \] 
then by taking the limite \[ \lim_{k \rightarrow  \infty } \frac{(\alpha k+\beta)+\alpha}{\alpha K+\beta} = 1 \]
Proof : We find the following limit : \[ \lim_{k \rightarrow \infty } \frac{\alpha +k \alpha \beta }{k \alpha + \beta } \] the leading term in the denominator of $ (\alpha + k \alpha \beta )/(k \alpha +\beta ) $ is k . Divide the numerator and denominator by this : \[  \lim_{ k \rightarrow \infty } \frac{\alpha/k + \alpha +\beta/k}{\alpha +\beta/k} \] 
the expressions $ \alpha/k $ and $ \beta/k $ both tend to zero as k approaches $\infty$ : $ \alpha/\alpha = 1 $ and therfore we conclude that  \[ \lim_{k \rightarrow  \infty } \frac{(\alpha k+\beta )+\alpha}{\alpha k+\beta } = 1 \] this means that the right-hand side also tends to 1 , So $ \lim_{k \rightarrow \infty } r_k=0 $ since the $ r_k  $  are bounded . We conclude  therefore the number $A(\alpha , \beta) $ is irrational since we know the   $ r_{k} > 1/b $ for all k . 
We use the above theorem to demonstrate the irrationality of several interesting constants. For the first case and by using Bertrand's postulate it is already shown that e is irrational by considering the sequence of integers 2, 3, 4, 5,... and the ratio of consecutive numbers must tend to one.

Setting $\alpha =2 $ and $\beta=1$ we get :
\[ \sum_{k=1}^{\infty} \frac{(2 k + 1) -1)} {\prod_{i=0}^{k-1} (2 k + 1) } =   \sqrt{2 \pi e } \hspace{1mm} \textbf{$erf$} \left(\frac{1}{\sqrt{2}} \right) +1 \] 
let $ O_{k}$ be the nth odd number , then by Bertrand's postulate we have $O_{k-1}  < g_{k} < O_{k}$ for all $k$, we may write $g_{k}=O_{k}-r_{k}$,where 0< $r_{k}$  < 1.Assume $g_{k}$ is rational such that \[ g_{k} =\frac{a}{b}\], where a,b are integers .Using the recurence relation 
\[ f_{k}=O_{k-1}\left(f_{k-1}-O_{k-1}+1\right) \]
it cleary seen that $b f_{k}$ is an integer for all k.
we write the above expression as follows
\[ \left(r_{k}+1\right)=\frac{f_{k+1}}{O_{k}}=\frac{O_{k+1}+r_{k+1}}{O_{K}} \]Since the $r_{k}$ are bounded, the right-hand side tends to 1 . But the  $ \lim _{k \rightarrow \infty} r_{k}=0 $  so we get the contradiction since $r_{k}  > 1/b $ and it easy to show that the limite  \[ \lim_{k \rightarrow  \infty } \frac{(2k+1)+2}{(2k+1)}=1\]  
More Exemples :
\[ A(3 , 1 ) = \sum_{k=1}^{\infty} \frac{(3 k + 1) -1)} {\prod_{i=0}^{k-1} (3 k + 1) } =   1+\sqrt[3]{3e} \Gamma\left(\frac{1}{3}\right)-\sqrt[3]{3e}  \Gamma \left(\frac{1}{3},\frac{1}{3}\right) \notin \mathbf{Q}\] 

\[ A(4,1) = \sum_{k=0}^{\infty} \frac{(4 k+2)-1}{\prod_{k=0}^{k-1}(4 k+2)}=2 \sqrt[4]{e} \sqrt{\pi} \hspace{1mm} {erf}\left(\frac{1}{2}\right)+2 \notin \mathbf{Q}\]
when $\alpha = i $ and $\beta =1 $  we get :
\[ A(i,1) = \sum_{k=0}^{\infty} \frac{(ki+1)-1)}{\prod_{i=0}^{k-1} (ki+1)} =
-\frac{i \Gamma(-i)}{\Gamma(1-i)}- \] \[ 
\frac{i\left(-\frac{i}{e}\right)^{i} \Gamma(-i)^{2}}{\Gamma(1-i)}+
\frac{i\left(-\frac{i}{e}\right)^{i} \Gamma(-i) \Gamma(-i,-i)}{\Gamma(1-i)}  \notin \mathbf{Q} \]

\[ A(5,3) = \sum_{k=0}^{\infty} \frac{(5k+3)-1}{\prod_{i=0}^{k-1} (5k+3)} 
1-\frac{4 \sqrt[5]{e} \Gamma\left(-\frac{2}{5}\right)}{5 \times 5^{2 / 5}}+ \\ \] \[ 
\frac{3 \sqrt[5]{e} \Gamma\left(\frac{3}{5}\right)}{5^{2 / 5}}+ \\
\frac{4 \sqrt[5]{e} \Gamma\left(-\frac{2}{5}, \frac{1}{5}\right)}{5 \times 5^{2 / 5}}- \\
\frac{3 \sqrt[5]{e} \Gamma\left(\frac{3}{5}, \frac{1}{5}\right)}{5^{2 / 5}}  \notin \mathbf{Q}\]

\bibliographystyle{unsrt}  
\bibliography{references}

Email:chtatbiali@gmail.com
\end{document}